\documentclass[11pt]{amsart}
\usepackage[margin=1.1in]{geometry}
\usepackage{amsmath,amssymb,amsthm}
\usepackage{booktabs}
\usepackage{xcolor}
\usepackage{tikz}
\usetikzlibrary{arrows.meta,positioning,calc}
\usepackage[colorlinks=true,linkcolor=blue,citecolor=blue,urlcolor=blue]{hyperref}

\makeatletter
\def\subsection{\@startsection{subsection}{2}%
  \z@{-.5\linespacing\@plus-.7\linespacing}{.5\linespacing}%
  {\normalfont\bfseries\centering}}
\makeatother
\newtheorem{theorem}{Theorem}[section]
\newtheorem{lemma}[theorem]{Lemma}
\newtheorem{proposition}[theorem]{Proposition}
\newtheorem{corollary}[theorem]{Corollary}
\newtheorem{conjecture}[theorem]{Conjecture}
\theoremstyle{definition}

\newtheorem{remark}[theorem]{Remark}
\newtheorem{example}[theorem]{Example}

\newcommand{\re}{\operatorname{re}}
\newcommand{\ro}{\operatorname{ro}}
\newcommand{\sgn}{\operatorname{sgn}}
\newcommand{\cyc}{\operatorname{cyc}}
\newcommand{\rad}{\operatorname{rad}}
\newcommand{\POW}{\mathrm{POW}}
\newcommand{\N}{\mathbb{N}}
\newcommand{\dd}[2]{(\!(#1,#2)\!)}

\begin{document}

\title{On the parity of $k$th roots of permutations}

\author{Luis Manuel Rivera}
\address{Unidad Acad\'emica de Matem\'aticas, Universidad Aut\'onoma de Zacatecas, Calzada Solidaridad y Paseo La Bufa, Col.~Hidr\'aulica, C.P.~98060, Zacatecas, Zac., Mexico}
\email{luismanuel.rivera@gmail.com}


\begin{abstract}
In every symmetric group $S_n$ with $n\ge2$, even and odd permutations are
equinumerous, but the $k$th roots of a given permutation rarely are. Let $k$ be even
and let $\sigma$ be a permutation with $k$th roots and $c_1$ fixed points. We show
that $\sigma$ has as many even as odd $k$th roots if and only if $T_k(c_1)=0$, where
$T_k(c)$ is the signed number of solutions of $x^k=1$ in $S_c$; in particular, the
cycles of length at least $2$ play no role. Congruences for $T_k$ modulo odd primes,
together with St\o rmer's theorem on consecutive smooth numbers, show that $T_k$ has
only finitely many zeros; we conjecture that $T_k(c)=0$ implies $c(c-1)\mid k$, and
prove this for $c\le40$. We also characterize the permutations whose $k$th roots all
have the same parity, and hence the $k$th powers in the alternating group, and derive
exact periodicity properties of the proportion of $m$th powers in $A_n$ when $m$ is a
prime power.
\end{abstract}
\maketitle

\section{Introduction}

Exactly half of the permutations of $[n]$, $n\ge2$, are even. Is the same true of the
square roots of a given permutation $\sigma$, that is, of the permutations $\tau$ with
$\tau^2=\sigma$? Almost never. The identity of $S_3$ has four square roots, namely the
identity and the three transpositions, and only one of them is even. On the other hand,
$(1)(2)(3\,4\,5)$ has exactly two square roots, $(3\,5\,4)$ and $(1\,2)(3\,5\,4)$, one of
each parity. What makes the second example balanced is its number of fixed points: as we
show, a permutation with square roots has as many even as odd square roots if and only if
it has exactly two fixed points.

The same phenomenon governs $k$th roots for every even $k$ (for $k$ odd, every $k$th root
of $\sigma$ has the parity of $\sigma$). The difference between the numbers of even and
odd $k$th roots of $\sigma$ factors over the cycle lengths of $\sigma$, and what had not
been observed is that the cycles of length at least $2$ never make it vanish: if $\ell$ is
even, all $k$th roots of a permutation made of $\ell$-cycles have the same sign, and if
$\ell\ge3$ is odd, the corresponding signed count is a positive multiple of an integer
congruent to $1$ modulo some odd $L\ge3$, and hence it is nonzero. Therefore the balance is decided by the fixed points alone, through a single
number: the signed count $T_k(c_1)$ of solutions of $x^k=1$ in the symmetric group on the
$c_1$ fixed points of $\sigma$ (Theorem~\ref{thm:balance}).

This raises the question of when $T_k(c)$ vanishes. It always vanishes at $c=2$, and if
$k$ is a power of $2$ only there; for other $k$ there are further zeros, such as
$T_6(3)=0$ and $T_{12}(4)=0$. Nothing in the definition of $T_k(c)$ as a signed sum
suggests that there should be only finitely many of them, and yet this is the case.
Congruences modulo the odd primes not dividing $k$ force every prime factor of $c(c-1)$
to divide $k$, so that $c-1$ and $c$ are consecutive integers whose prime factors lie in
a fixed finite set; by a classical theorem of St\o rmer there are only finitely many such
pairs (Corollary~\ref{cor:radical}). We conjecture that, in fact, $c(c-1)$ divides $k$
(Conjecture~\ref{conj:T}), and we prove this for $c\le40$
(Proposition~\ref{prop:conj40}).

Roots of permutations have a long history. The permutations having $k$th roots are
characterized in \cite[Theorem~4.43, p.~158]{Wilf} (see also \cite{AJS,BC}), the number of
$k$th roots is computed in \cite{LMR,Pavlov}, and the probability that a random
permutation has a $k$th root is studied in \cite{BMW,CW,Chernoff94,LMR,Pouyanne}. On
parity, Moser and Wyman \cite{MW} counted the even solutions of $x^2=1$, Chernoff
\cite{Chernoff90} enumerated by cycle type the solutions of $x^r=\alpha$ in $A_n$,
Pournaki \cite{Pournaki} counted the squares of $A_n$, and generating functions for the
numbers $\re_k(\sigma)$ and $\ro_k(\sigma)$ of even and odd $k$th roots were obtained in
\cite{GLR}. Adin, Heged\H{u}s and Roichman \cite{AHR} studied the \emph{signed root
enumerator}
\[
\bar\rho_k(\sigma)=\sum_{\tau^k=\sigma}\sgn(\tau)=\re_k(\sigma)-\ro_k(\sigma),
\]
a virtual character of $S_n$, and its analogues for the Weyl groups of types $B$ and $D$
\cite{AHR2}; and Singh \cite{Singh} related the odd square roots of even permutations to
the sum of element orders in $S_n$ and $A_n$. Results of \cite{Chernoff90}, \cite{Pournaki}
and \cite{Singh} are recovered below as special cases.

\subsection{Notation}

For $\sigma\in S_n$ we write $c_\ell=c_\ell(\sigma)$ for the number of $\ell$-cycles of $\sigma$. Since an $\ell$-cycle is a product of $\ell-1$ transpositions, $\sgn(\sigma)=\prod_\ell(-1)^{(\ell+1)c_\ell}$; in particular, cycles of odd length are even permutations, and $\sgn(\sigma)=(-1)^{\#\{\text{cycles of even length}\}}$. We say that $\sigma$ has \emph{type} $(\ell)^c$ if all its cycles have length $\ell$ and there are $c$ of them. For $\ell,k\in\N$ let
\[
G_k(\ell)=\{g\in\N:\gcd(g\ell,k)=g\},\qquad \dd{\ell}{k}=\prod_{p\mid\ell}p^{\nu_p(k)},
\]
where $\nu_p$ is the $p$-adic valuation and $\dd{1}{k}=1$. Recall \cite[Section~2]{LMR} that $\dd{\ell}{k}$ is the least element of $G_k(\ell)$ (see also Lemma~\ref{lem:structure} below), and that $\sigma$ has a $k$th root if and only if $\dd{\ell}{k}\mid c_\ell$ for every $\ell$ \cite[Theorem~4.43, p.~158]{Wilf}. Every $k$th root of $\sigma$ is obtained by grouping, for each $\ell$ and each $g\in G_k(\ell)$, disjoint sets of $g$ cycles of length $\ell$ of $\sigma$ into cycles of length $g\ell$ \cite[Section~3]{LMR}; see Lemma~\ref{lem:cyclepower} and Figure~\ref{fig:cycle}.

In this note we study when roots of both parities occur and when they are balanced. If $k$ is odd, then $\sgn(\tau^k)=\sgn(\tau)^k=\sgn(\tau)$, so every $k$th root of $\sigma$ has the parity of $\sigma$; hence we assume that $k$ is even. For $c\ge 0$ let
\[
T_k(c)=\sum_{\substack{\pi\in S_c\\ \pi^k=1}}\sgn(\pi)
\]
be the signed number of solutions of $x^k=1$ in $S_c$ (with $S_0$ the trivial group, so $T_k(0)=1$). For $m\ge1$ let $p_m(n)$ be the probability that a random element of $S_n$ has an $m$th root, and for $n\ge2$ let $p^A_m(n)=|A_n^m|/|A_n|$, where $A_n^m=\{\tau^m:\tau\in A_n\}$. For $m\ge1$ write
 \[
\exp_m(y)=\sum_{i\ge0}\frac{y^{im}}{(im)!}
\]
and
 \[
 \widetilde{\exp}_m(y)=\sum_{i\ge0}\frac{(-1)^iy^{im}}{(im)!}.
 \]
Thus $\exp_1(y)=e^y$, $\exp_2(y)=\cosh(y)$ and $\widetilde{\exp}_2(y)=\cos(y)$.

\subsection{Main results}

Our first result says that, for $k$ even, the balance between even and odd roots is decided by the fixed points alone.

\begin{theorem}\label{thm:balance}
Let $k$ be even and let $\sigma$ be a permutation with $k$th roots and $c_1$ fixed points. Then
\[
\re_k(\sigma)=\ro_k(\sigma)\iff T_k(c_1)=0 .
\]
In particular $\re_k(\sigma)\neq\ro_k(\sigma)$ if $c_1\le1$, and $\re_k(\sigma)=\ro_k(\sigma)$ if $c_1=2$. Equivalently, $\bar\rho_k(\sigma)=0$ if and only if $\sigma$ has no $k$th root or $T_k(c_1(\sigma))=0$.
\end{theorem}

Our second result concerns congruences for $T_k$. Let $r_k(c)$ and $e_k(c)$ be the numbers of solutions of $x^k=1$ in $S_c$ and in $A_c$, respectively, so that $e_k(c)=\tfrac12(r_k(c)+T_k(c))$. These congruences are the key to the zeros of $T_k$.

\begin{theorem}\label{thm:congr}
Let $p$ be an odd prime and $c\ge p$.
\begin{enumerate}
\item If $p\nmid k$, then, modulo $p$,
\begin{enumerate}
\item $r_k(c)\equiv r_k(c-p)$,
\item $T_k(c)\equiv T_k(c-p)$,
\item $e_k(c)\equiv e_k(c-p)$.
\end{enumerate}
\item If $p\mid k$, then $r_k(c)\equiv T_k(c)\equiv e_k(c)\equiv0\pmod p$.
\end{enumerate}
\end{theorem}

For $c=p$ we recover the result of Chernoff \cite{Chernoff90}: the number of solutions of $x^k=1$ in $A_p$ is congruent modulo $p$ to $1$ if $p\nmid k$ and to $0$ if $p\mid k$ (as noted in \cite[p.~227]{Chernoff90}, the first of these congruences appears in \cite{MW} with the condition $k\nmid p$ instead of $p\nmid k$).

\begin{corollary}\label{cor:radical}
Let $k$ be even and $c\ge2$. If $T_k(c)=0$, then every prime divisor of $c(c-1)$ divides $k$. Consequently, for fixed $k$ the set $\{c: T_k(c)=0\}$ is finite, and if $k$ is a power of $2$ then $T_k(c)=0$ if and only if $c=2$.
\end{corollary}
The finiteness of the zeros of $T_k$ is not apparent from its definition as a signed
sum; it follows from the congruences of Theorem~\ref{thm:congr} combined with
St\o rmer's theorem on consecutive smooth numbers, as we show in Section~\ref{sec:congr}.
\begin{conjecture}\label{conj:T}
If $k$ is even, $c\ge2$ and $T_k(c)=0$, then $c(c-1)\mid k$.
\end{conjecture}
Corollary~\ref{cor:radical} proves Conjecture~\ref{conj:T} at the level of radicals:
it shows that $\rad(c(c-1))\mid k$ whenever $c\ge2$ and $T_k(c)=0$, where $\rad(n)$
denotes the product of the distinct prime divisors of $n$. Moreover, for a fixed set
of primes dividing $k$, it reduces the determination of the zeros of $T_k$ to a finite
computation. For example, if $k=2^a3^b$ with $a\ge1$, then $T_k(c)=0$ if and only if $c=2$, or $c=3$ and $6\mid k$, or $c=4$ and $12\mid k$, or $c=9$ and $72\mid k$ (see Table~\ref{tab:23}). In the other direction, for each fixed $c$ the conjecture reduces to a finite computation, which we have carried out for $c\le40$.

\begin{proposition}\label{prop:conj40}
Conjecture~\ref{conj:T} holds for $2\le c\le40$: if $k$ is even, $2\le c\le 40$ and
$T_k(c)=0$, then $c(c-1)\mid k$.
\end{proposition}

The reason is that, for fixed $c$, both $T_k(c)$ and the condition $c(c-1)\mid k$ depend
only on $\gcd(k,\operatorname{lcm}(1,\dots,c))$; see Section~\ref{sec:congr}, where
Table~\ref{tab:zeros} also lists, for $3\le c\le40$, the least even $k$ with $T_k(c)=0$.
The condition $c(c-1)\mid k$ is not sufficient, since $T_{56}(8)=2688\neq0$.

Our next result is a characterization of the permutations all of whose $k$th roots
have the same parity. Let $a=\nu_2(k)$. If $\sigma$ has $k$th roots and $\ell$ is even, then $2^a$ divides $\dd{\ell}{k}$
(since $2\mid\ell$), and $\dd{\ell}{k}$ divides $c_\ell(\sigma)$; hence
$c_\ell(\sigma)/2^a$ is an integer.

\begin{theorem}\label{thm:rigid}
Let $k$ be even, $a=\nu_2(k)$, and let $\sigma$ be a permutation with $k$th roots.
All $k$th roots of $\sigma$ have the same parity if and only if
$c_\ell(\sigma)\in\{0,\dd{\ell}{k}\}$ for every odd $\ell$. In that case, they all
have sign $(-1)^{\sum_{\ell\ \mathrm{even}}c_\ell(\sigma)/2^a}$.
\end{theorem}

A permutation $\sigma\in A_n$ is a $k$th power in $A_n$ if $\sigma=\tau^k$ for some $\tau\in A_n$. If $k$ is odd, $\sgn(\tau^k)=\sgn(\tau)$, so the $k$th powers in $A_n$ are exactly the even $k$th powers in $S_n$. For $k$ even every $k$th power is even, and Theorem~\ref{thm:rigid} gives:

\begin{corollary}\label{cor:An}
Let $k$ be even and $a=\nu_2(k)$. A permutation $\sigma$ is the $k$th power of an even permutation if and only if $\dd{\ell}{k}\mid c_\ell(\sigma)$ for every $\ell$, and either $c_\ell(\sigma)>\dd{\ell}{k}$ for some odd $\ell$, or $\sum_{\ell\ \mathrm{even}}c_\ell(\sigma)/2^a$ is even.
\end{corollary}

For $k=2$, this recovers Pournaki's characterization of the squares in $A_n$ \cite{Pournaki}. Theorem~\ref{thm:rigid} also contains \cite[Lemmas~3.1 and 3.2]{Singh}: a permutation with square roots and $c_1\ge2$ has square roots of both parities, and an even permutation with pairwise distinct cycle lengths has no odd square roots.

The reduction behind Theorems~\ref{thm:balance} and~\ref{thm:rigid} also
characterizes the permutations with a unique $k$th root (Corollary~\ref{cor:unique} in
Section~\ref{sec:rigid}); for $k$ even they are an extreme case of
Theorem~\ref{thm:rigid}.

Corollary~\ref{cor:An}, together with the identity $\sgn(\tau^m)=\sgn(\tau)$ for $m$ odd,
describes the set $A_n^m$ of $m$th powers in $A_n$ explicitly. This allows us to compute
the proportion $p^A_m(n)$ exactly when $m$ is a prime power, and hence to compare it with
the following periodicity of $p_m(n)$. It was proved in~\cite{Chernoff94,LMR} (see \cite{BMW} for $m$ prime) that if $m=q^r$ with $q$ prime, then
\[
p_m(qj)=p_m(qj+1)=\cdots=p_m(qj+q-1) \quad \text{ for } j\ge0.
\]
For $A_n$ we prove the following analogues; the second one explains a failure of periodicity observed by Pournaki.

\begin{theorem}\label{thm:Aodd}
Let $m=p^r$ with $p$ an odd prime. Define $\kappa_j$ by
\[
\sum_{j\ge0}\kappa_j\,x^{pj}=(1+x^p)^{-1/p}\prod_{i\ge1}\exp_m\!\left(\frac{(-1)^{i+1}x^{pi}}{pi}\right).
\]
(so $\kappa_0=1$). Then, for every $j\ge0$ and $n=pj+i\ge2$ with $0\le i\le p-1$,
\[
p^A_m(n)=\begin{cases}p_m(pj)+\kappa_j,& i\in\{0,1\},\\ p_m(pj),& 2\le i\le p-1.\end{cases}
\]
\end{theorem}

\begin{theorem}\label{thm:A2r}
Let $m=2^r$ with $r\ge1$. Then, for $n\ge2$,
\[
p^A_m(n)=2p_m(n)-2\,[x^n]\,O(x)\,\frac{E(x)-\widetilde E(x)}{2},
\]
where
\begin{align*}
O(x)&=\prod_{\ell\ \mathrm{odd}}(1+x^\ell/\ell),\\
 E(x)&=\prod_{\ell\ \mathrm{even}}\exp_m(x^\ell/\ell),\\
 \widetilde E(x)&=\prod_{\ell\ \mathrm{even}}\widetilde{\exp}_m(x^\ell/\ell).
 \end{align*}
 Moreover $p^A_m(2j)=p^A_m(2j+1)$ for $1\le j\le m$, and
\[
p^A_m(2m+2)-p^A_m(2m+3)=\frac{2}{3\cdot 2^m\,m!}.
\]
\end{theorem}

For $m=2$, Theorem~\ref{thm:A2r} explains the failure $p^A_2(6)=3/4\ne p^A_2(7)=2/3$ observed by Pournaki \cite{Pournaki}. Moreover, in contrast with the monotonicity of $p_m(n)$ for prime powers ($m$ prime \cite{BMW}, prime powers \cite{CW}), $p^A_m(n)$ need not be monotone: $p^A_3(4)=1/3<p^A_3(5)=2/3$, $p^A_5(6)=3/5<p^A_5(7)=4/5$ and $p^A_2(8)=9/16<p^A_2(9)=143/240$.

Section~\ref{sec:prelim} contains the preliminaries; the remaining sections contain
the proofs, in the order in which the results are stated.

\section{Preliminaries}\label{sec:prelim}

\subsection{Powers of cycles and the structure of roots}
We use the following well-known facts (see, e.g., \cite[Section~3]{LMR}).

\begin{lemma}\label{lem:cyclepower}
Let $C$ be a cycle of length $N$ and let $k\ge1$. Then $C^k$ is a product of $\gcd(N,k)$ disjoint cycles, each of length $N/\gcd(N,k)$. Consequently, if $\tau^k=\sigma$, then each cycle of $\tau$ of length $g\ell$ with $\gcd(g\ell,k)=g$ splits into exactly $g$ cycles of $\sigma$ of length $\ell$, and conversely every cycle of $\tau$ whose $k$th power consists of $\ell$-cycles has length $g\ell$ for some $g\in G_k(\ell)$.
\end{lemma}

\begin{figure}[ht]
\centering
\begin{tikzpicture}[scale=1.0,
  pt/.style={circle, fill=black, inner sep=1.4pt},
  arr/.style={-{Stealth[length=2mm]}, thick}]
\begin{scope}
\foreach \i in {1,...,6} {
  \node[pt] (a\i) at ({90-60*(\i-1)}:1.3) {};
  \node at ({90-60*(\i-1)}:1.65) {$\i$};
}
\foreach \i/\j in {1/2,2/3,3/4,4/5,5/6,6/1} \draw[arr, blue] (a\i) to[bend left=18] (a\j);
\node at (0,-2.2) {$\tau=(1\,2\,3\,4\,5\,6)$};
\end{scope}
\draw[-{Stealth[length=3mm]}, very thick] (2.2,0) -- node[above] {$x\mapsto x^2$} (3.6,0);
\begin{scope}[xshift=5.8cm]
\foreach \i in {1,...,6} {
  \node[pt] (b\i) at ({90-60*(\i-1)}:1.3) {};
  \node at ({90-60*(\i-1)}:1.65) {$\i$};
}
\foreach \i/\j in {1/3,3/5,5/1} \draw[arr, red!80!black] (b\i) -- (b\j);
\foreach \i/\j in {2/4,4/6,6/2} \draw[arr, green!50!black, dashed] (b\i) -- (b\j);
\node at (0,-2.2) {$\tau^2=(1\,3\,5)(2\,4\,6)$};
\end{scope}
\end{tikzpicture}
\caption{Lemma~\ref{lem:cyclepower} with $N=6$, $k=2$: here $\ell=3$, $g=2\in G_2(3)$, and the $6$-cycle $\tau$ splits into $g=2$ cycles of length $\ell=3$. Conversely, a square root of $(1\,3\,5)(2\,4\,6)$ is obtained by \emph{interleaving} the two $3$-cycles; there are $(g-1)!\,\ell^{g-1}=3$ ways to do so, namely $(1\,2\,3\,4\,5\,6)$, $(1\,4\,3\,6\,5\,2)$ and $(1\,6\,3\,2\,5\,4)$.}
\label{fig:cycle}
\end{figure}
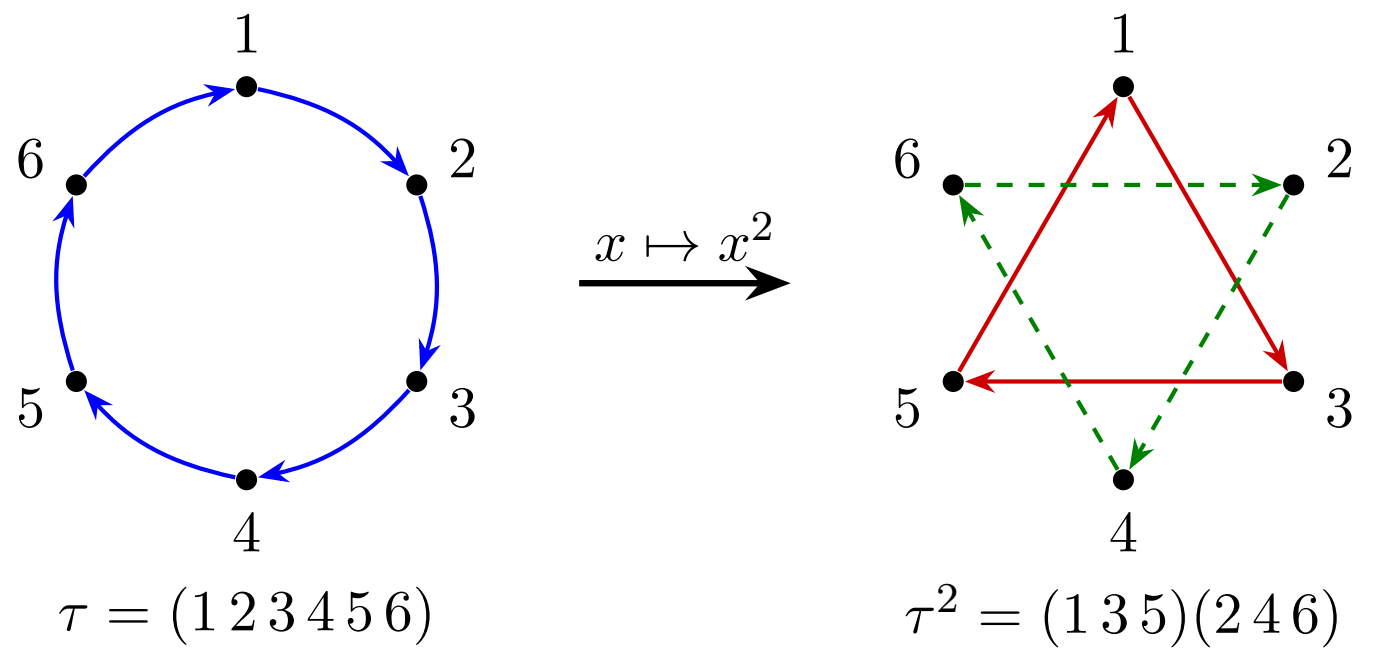
For $c\ge0$ let $r_k(\ell,c)$, $\re_k(\ell,c)$ and $\ro_k(\ell,c)$ be the numbers of $k$th
roots, even $k$th roots and odd $k$th roots of a fixed permutation of type $(\ell)^c$ (these
numbers do not depend on the chosen permutation, since conjugation preserves the sign), and
\[
D_k(\ell,c):=\re_k(\ell,c)-\ro_k(\ell,c).
\]
For every $\sigma$ we have (see \cite[Section~3]{LMR} and \cite[Lemma~2.5]{GLR})
\begin{equation}\label{eq:product}
r_k(\sigma)=\prod_{\ell}r_k(\ell,c_\ell(\sigma)),\qquad
\bar\rho_k(\sigma)=\re_k(\sigma)-\ro_k(\sigma)=\prod_{\ell}D_k(\ell,c_\ell(\sigma)),
\end{equation}
where $r_k(\sigma)$ is the number of $k$th roots of $\sigma$. Clearly $|D_k(\ell,c)|\le r_k(\ell,c)$,
with equality if and only if all $k$th roots of a permutation of type $(\ell)^c$ have the same
parity.

\subsection{The exponential formula}

We use the following form of the exponential formula. It follows by combining \cite[Corollary~5.1.8]{Stanley}, applied separately to the cycles of each length $\ell$ (taking there $f$ to be the indicator function of $\{\ell\}$, and the second function to be $i\mapsto a_\ell(i)$), with the product formula \cite[Proposition~5.1.3]{Stanley}; the case $a_\ell(i)=w_\ell^{\,i}$ for all $\ell$ and $i$ (no restrictions on the numbers of cycles)
is \cite[Corollary~5.1.9]{Stanley}, with $f(n)=w_n$. See also \cite[Theorem~4.38, p.~152]{Wilf} and the argument leading to \cite[(4.73), p.~157]{Wilf}. For the convenience of the reader, we include the short direct proof.

\begin{lemma}\label{lem:expformula}
For each $\ell\ge1$ let $(a_\ell(i))_{i\ge0}$ be a sequence in a commutative ring containing
$\mathbb Q$, with $a_\ell(0)=1$. Then
\[
\sum_{c\ge0}\frac{t^c}{c!}\sum_{\pi\in S_c}\ \prod_{\ell\ge1}a_\ell\big(c_\ell(\pi)\big)
=\prod_{\ell\ge1}\ \sum_{i\ge0}\frac{a_\ell(i)}{i!}\Big(\frac{t^\ell}{\ell}\Big)^{i}.
\]
In particular, if $a_\ell(i)=w_\ell^{\,i}$ when $i\in A_\ell$ and $a_\ell(i)=0$ otherwise, where
$0\in A_\ell\subseteq\{0,1,2,\dots\}$, the left-hand side counts the permutations with
$c_\ell(\pi)\in A_\ell$ for all $\ell$, each weighted by $\prod_C w_{|C|}$; and if $A_\ell=\{0,1,2,\dots\}$
for all $\ell$, the right-hand side is $\exp\big(\sum_\ell w_\ell t^\ell/\ell\big)$.
\end{lemma}

\begin{proof}
There are $c!/\prod_\ell \ell^{c_\ell}c_\ell!$ permutations of $[c]$ with $c_\ell$ cycles of length $\ell$ for
each $\ell$. Hence the coefficient of $t^c/c!$ on the left-hand side is
$c!\sum\prod_\ell a_\ell(c_\ell)/(\ell^{c_\ell}c_\ell!)$, the sum over all $(c_\ell)_{\ell\ge1}$ with $\sum_\ell \ell c_\ell=c$,
which is the coefficient of $t^c/c!$ on the right-hand side.
\end{proof}

\subsection{\texorpdfstring{Generating functions for $r_k(\ell,c)$ and $D_k(\ell,c)$}{Generating functions for r\_k(l,c) and D\_k(l,c)}}
The following generating functions are known. The first one is the specialization to
permutations of type $(\ell)^c$ of the multivariate generating function in
\cite[Theorem~2]{LMR} (put $t_\ell=t$ and $t_j=0$ for $j\ne\ell$), and the second one is
\cite[Lemma~2.10]{GLR}, with $t_\ell$ written as $t$.

\begin{proposition}[{\cite[Theorem~2]{LMR}, \cite[Lemma~2.10]{GLR}}]\label{prop:egf}
For all $k,\ell\ge1$,
\begin{align}
\sum_{c\ge0}r_k(\ell,c)\frac{t^c}{c!}&=\exp\Big(\sum_{g\in G_k(\ell)}\frac{\ell^{g-1}}{g}t^g\Big),\label{eq:egf}\\
\sum_{c\ge0}D_k(\ell,c)\frac{t^c}{c!}&=\exp\Big(\sum_{g\in G_k(\ell)}(-1)^{g\ell+1}\frac{\ell^{g-1}}{g}t^g\Big).\label{eq:egfD}
\end{align}
\end{proposition}

We also need the following count \cite[Proposition~7]{LMR}: if $g\in G_k(\ell)$, then every permutation of type $(\ell)^g$ is the $k$th power of exactly $(g-1)!\,\ell^{g-1}$ cycles of length $g\ell$ (see Figure~\ref{fig:cycle} for $g=2$, $\ell=3$). Consequently, if moreover $g\mid c$, the number of $k$th roots of a permutation of type $(\ell)^c$ all of whose cycles have length $g\ell$ (choose a partition of the $c$ cycles into blocks of size $g$, and then one $g\ell$-cycle for each block) is
\begin{equation}\label{eq:onetype}
\frac{c!}{(g!)^{c/g}(c/g)!}\big((g-1)!\,\ell^{g-1}\big)^{c/g}=\frac{c!\,\ell^{(g-1)c/g}}{g^{c/g}\,(c/g)!}.
\end{equation}
Note that $c!/(g^{c/g}(c/g)!)$ is a positive integer: it is the number of permutations of $[c]$ of type $(g)^{c/g}$.

\subsection{Three technical lemmas}
The next three lemmas reduce the computation of $r_k(\ell,c)$ and $D_k(\ell,c)$ to the
case in which $\ell$ is coprime to $k$, and relate the latter to the solutions of $x^k=1$
in a symmetric group; in particular, they explain how $T_k$ enters the paper.

\begin{lemma}\label{lem:structure}
Let $g_0=\dd{\ell}{k}$ and $k'=k/g_0$. Then $\gcd(k',\ell g_0)=1$ and $G_k(\ell)=\{g_0e: e\mid k'\}$. In particular, $g_0$ is the least element of $G_k(\ell)$ and every element of $G_k(\ell)$ is a multiple of $g_0$.
\end{lemma}

\begin{proof}
If $g\in G_k(\ell)$, then $g=\gcd(g\ell,k)$ divides $k$. Conversely, for $g\mid k$ we have
\[
\gcd(g\ell,k)=\gcd\big(g\ell,\,g\cdot(k/g)\big)=g\cdot\gcd(\ell,k/g).
\]
Hence
\begin{equation}\label{eq:Gk}
G_k(\ell)=\{g:\ g\mid k\ \text{ and }\ \gcd(\ell,k/g)=1\}
\end{equation}
(this is \cite[Remark~1]{LMR}). Write $k=g_0k'$ with $g_0=\prod_{p\mid\ell}p^{\nu_p(k)}$ and $k'=\prod_{p\nmid\ell}p^{\nu_p(k)}$. A divisor $d$ of $k$ satisfies $\gcd(\ell,d)=1$ if and only if no prime divisor of $\ell$ divides $d$, that is, if and only if $d\mid k'$. Applying this to $d=k/g$, we get $g\in G_k(\ell)$ if and only if $k/g\mid k'$, i.e.\ if and only if $g=g_0e$ with $e=k'/(k/g)$ a divisor of $k'$. Finally, the prime divisors of $\ell g_0$ are exactly those of $\ell$, and none of them divides $k'$; so $\gcd(k',\ell g_0)=1$.
\end{proof}

For instance, if $k=12$ and $\ell=3$, then $g_0=3$, $k'=4$ and $G_{12}(3)=\{3,6,12\}$: a $3$-cycle of $\sigma$ can only come from cycles of length $9$, $18$ or $36$ of a $12$th root.

\begin{lemma}\label{lem:reduction}
Keep the notation of Lemma~\ref{lem:structure} and let $L=\ell g_0$. Then $r_k(\ell,c)=0$ unless $g_0\mid c$, and for $j\ge0$
\[
r_k(\ell,g_0j)=\frac{(g_0j)!}{j!}\,\frac{\ell^{(g_0-1)j}}{g_0^{\,j}}\;r_{k'}(L,j).
\]
If $\ell$ is odd, the same identity holds with $r$ replaced by $D$.
\end{lemma}

\begin{proof}
By Lemma~\ref{lem:structure}, the exponent in \eqref{eq:egf} is
\[
\sum_{g\in G_k(\ell)}\frac{\ell^{g-1}}{g}t^g=\sum_{e\mid k'}\frac{\ell^{g_0e-1}}{g_0e}\,t^{g_0e}.
\]
Put $v=(\ell t)^{g_0}/L$. Then, for every $e\ge1$,
\[
\frac{L^{e-1}v^e}{e}=\frac{(\ell g_0)^{e-1}\,\ell^{g_0e}\,t^{g_0e}}{(\ell g_0)^{e}\,e}=\frac{\ell^{g_0e-1}}{g_0e}\,t^{g_0e}.
\]
Since $\gcd(L,k')=1$ by Lemma~\ref{lem:structure}, for every $e\mid k'$ we have $\gcd(eL,k')=e\gcd(L,k'/e)=e$; by \eqref{eq:Gk} (applied to $k'$ and $L$), $G_{k'}(L)=\{e:e\mid k'\}$. Therefore, using \eqref{eq:egf} twice,
\[
\sum_{c\ge0}r_k(\ell,c)\frac{t^c}{c!}=\exp\Big(\sum_{e\in G_{k'}(L)}\frac{L^{e-1}}{e}v^e\Big)=\sum_{j\ge0}r_{k'}(L,j)\frac{v^j}{j!}=\sum_{j\ge0}r_{k'}(L,j)\frac{\ell^{g_0j}}{L^j\,j!}\,t^{g_0j}.
\]
The right-hand side contains only powers $t^c$ with $g_0\mid c$, so $r_k(\ell,c)=0$ if $g_0\nmid c$. Comparing the coefficients of $t^{g_0j}$ gives
\[
r_k(\ell,g_0j)=(g_0j)!\,\frac{\ell^{g_0j}}{(\ell g_0)^j\,j!}\,r_{k'}(L,j)=\frac{(g_0j)!}{j!}\,\frac{\ell^{(g_0-1)j}}{g_0^{\,j}}\,r_{k'}(L,j).
\]
If $\ell$ is odd, then every prime divisor of $g_0$ divides $\ell$, so $g_0$ and $L$ are odd. Hence, for $g=g_0e$, $(-1)^{\ell g+1}=(-1)^{\ell g_0e+1}=(-1)^{e+1}=(-1)^{Le+1}$, so the sign factors in \eqref{eq:egfD} for $(k,\ell)$ and for $(k',L)$ correspond under the same substitution, and the computation above applies verbatim to $D$.
\end{proof}

\begin{lemma}\label{lem:coprime}
If $\gcd(L,k)=1$, then
\[
r_k(L,c)=\sum_{\substack{\pi\in S_c\\ \pi^k=1}}L^{\,c-\cyc(\pi)},\qquad\text{and, if $L$ is odd,}\qquad D_k(L,c)=\sum_{\substack{\pi\in S_c\\ \pi^k=1}}\sgn(\pi)\,L^{\,c-\cyc(\pi)},
\]
where $\cyc(\pi)$ is the number of cycles of $\pi$. In particular $D_k(1,c)=T_k(c)$ and, for $L$ odd, $D_k(L,c)\equiv1\pmod L$.
\end{lemma}

\begin{proof}
By \eqref{eq:Gk}, $G_k(L)=\{g:g\mid k,\ \gcd(L,k/g)=1\}$ is the set of all divisors of $k$, because $\gcd(L,k)=1$. Take $A_g=\{0,1,2,\dots\}$ if $g\mid k$, $A_g=\{0\}$ otherwise, in Lemma~\ref{lem:expformula}; the permutations counted are exactly those whose cycle lengths divide $k$, that is, the $\pi\in S_c$ with $\pi^k=1$. We obtain
\[
c!\,[t^c]\exp\Big(\sum_{g\mid k}\frac{w_g}{g} t^g\Big)=\sum_{\substack{\pi\in S_c\\ \pi^k=1}}\ \prod_{C}w_{|C|}.
\]
Taking $w_g=L^{g-1}$ and comparing with \eqref{eq:egf} gives the formula for $r_k(L,c)$, because $\prod_C L^{|C|-1}=L^{\sum_C(|C|-1)}=L^{c-\cyc(\pi)}$. For $D$ and $L$ odd we take, according to \eqref{eq:egfD},
\[
w_g=(-1)^{Lg+1}L^{g-1}=(-1)^{g+1}L^{g-1},
\]
and use $\prod_C(-1)^{|C|+1}=\sgn(\pi)$. For $L=1$ we get $D_k(1,c)=\sum_{\pi^k=1}\sgn(\pi)=T_k(c)$. Finally, $\cyc(\pi)\le c$ with equality only for the identity, so every term with $\pi\neq1$ is divisible by $L$, while the identity contributes $\sgn(1)L^0=1$; hence $D_k(L,c)\equiv1\pmod L$.
\end{proof}

\section{Proof of Theorem~\ref{thm:balance}}\label{sec:balance}

The key observation is that, when $k$ is even, the cycles of even length of $\sigma$ do not influence the balance.

\begin{lemma}\label{lem:even-ell}
Let $k$ be even, $a=\nu_2(k)$, and $\ell$ even. Then every $k$th root of a permutation of type $(\ell)^c$ has sign $(-1)^{c/2^a}$. In particular $D_k(\ell,c)=(-1)^{c/2^a}r_k(\ell,c)$.
\end{lemma}

\begin{proof}
If $r_k(\ell,c)=0$ there is nothing to prove, so assume that a root exists.

\emph{Step 1: every $g\in G_k(\ell)$ satisfies $\nu_2(g)=a$.} Let $g\in G_k(\ell)$. Since $g\mid k$, $\nu_2(g)\le a$. Suppose $\nu_2(g)<a$. Since $2\mid\ell$,
\[
\nu_2(\gcd(g\ell,k))=\min\big(\nu_2(g)+\nu_2(\ell),\,a\big)\ge\nu_2(g)+1,
\]
contradicting $\gcd(g\ell,k)=g$. Hence $\nu_2(g)=a$, i.e.\ $g/2^a$ is odd.

\emph{Step 2: the sign of a root.} By Lemma~\ref{lem:cyclepower}, a root $\tau$ has cycle type $\prod_g(g\ell)^{s_g}$ for some non-negative integers $s_g$, $g\in G_k(\ell)$, with $\sum_g g s_g=c$ (each cycle of length $g\ell$ uses $g$ of the $c$ cycles of $\sigma$). All cycles of $\tau$ have even length $g\ell$, so $\sgn(\tau)=(-1)^{\sum_g s_g}$. By Step~1, $2^a\mid c$ and
\[
\frac{c}{2^a}=\sum_g\frac{g}{2^a}\,s_g\equiv\sum_g s_g\pmod 2,
\]
because each $g/2^a$ is odd. Thus $\sgn(\tau)=(-1)^{c/2^a}$ for every root $\tau$.
\end{proof}

\begin{example}
Let $k=2$ (so $a=1$) and $\sigma=(1\,2)(3\,4)$, of type $(2)^2$. Here $G_2(2)=\{2\}$, and the square roots of $\sigma$ are the $4$-cycles $(1\,3\,2\,4)$ and $(1\,4\,2\,3)$; both are odd, in accordance with $(-1)^{c/2^a}=(-1)^{1}$. In particular, $\sigma$ is not a square in $A_4$.
\end{example}

\begin{proof}[Proof of Theorem~\ref{thm:balance}]
By \eqref{eq:product}, $\bar\rho_k(\sigma)=\prod_\ell D_k(\ell,c_\ell)$ with $c_\ell=c_\ell(\sigma)$. Since $\sigma$ has $k$th roots, $r_k(\ell,c_\ell)>0$ for every $\ell$. We show that $D_k(\ell,c_\ell)\ne0$ for every $\ell\ge2$; then $\bar\rho_k(\sigma)=0$ if and only if $D_k(1,c_1)=0$, and $D_k(1,c_1)=T_k(c_1)$ by Lemma~\ref{lem:coprime} (with $L=1$; note that $\dd{1}{k}=1$).

\emph{Case $\ell$ even.} By Lemma~\ref{lem:even-ell}, $D_k(\ell,c_\ell)=\pm r_k(\ell,c_\ell)\ne0$.

\emph{Case $\ell\ge3$ odd.} Let $g_0=\dd{\ell}{k}$, $k'=k/g_0$ and $L=\ell g_0$. As $\sigma$ has roots, $g_0\mid c_\ell$ (Lemma~\ref{lem:reduction}); put $j=c_\ell/g_0$. By Lemma~\ref{lem:reduction},
\[
D_k(\ell,c_\ell)=\lambda\, D_{k'}(L,j),\qquad \lambda=\frac{(g_0j)!}{j!}\,\frac{\ell^{(g_0-1)j}}{g_0^{\,j}}>0.
\]
Now $L=\ell g_0\ge3$ is odd and $\gcd(L,k')=1$ (Lemma~\ref{lem:structure}), so Lemma~\ref{lem:coprime} gives $D_{k'}(L,j)\equiv1\pmod L$. Since $L\ge3$, this forces $D_{k'}(L,j)\ne0$, and hence $D_k(\ell,c_\ell)\neq0$.

This proves the equivalence. For the particular cases: $S_0$ and $S_1$ are trivial, so $T_k(0)=T_k(1)=1\ne0$; and $S_2=\{1,(1\,2)\}$, where both elements satisfy $x^k=1$ because $k$ is even, so $T_k(2)=1-1=0$.
\end{proof}

\begin{example}
Let $k=2$ and $\sigma=(1)(2)(3\,4\,5)$, so $c_1=2$ and $c_3=1$. The square roots of $\sigma$ are $(3\,5\,4)$ and $(1\,2)(3\,5\,4)$, one even and one odd, as predicted by Theorem~\ref{thm:balance}. In the notation of the proof, $D_2(3,1)=1$ and $D_2(1,2)=T_2(2)=0$.
\end{example}

\begin{remark}
For $c_1=2$ there is a direct bijective proof. Let $a,b$ be the fixed points of $\sigma$.
If $\tau^k=\sigma$, then by Lemma~\ref{lem:cyclepower} the support of the cycle of $\tau$
through $a$ consists of fixed points of $\sigma$, and similarly for $b$; hence $\tau$ either
fixes $a$ and $b$ or swaps them. In both cases $\tau$ commutes with $(a\ b)$, so
$(\tau\,(a\ b))^k=\tau^k(a\ b)^k=\sigma$ since $k$ is even. Thus $\tau\mapsto\tau\,(a\ b)$
is a parity-reversing involution on the set of $k$th roots of $\sigma$.
\end{remark}

\section{\texorpdfstring{Proof of Theorem~\ref{thm:congr}, Corollary~\ref{cor:radical} and Proposition~\ref{prop:conj40}}{Proof of Theorem 1.2, Corollary 1.3 and Proposition 1.5}}\label{sec:congr}

\begin{proof}[Proof of Theorem~\ref{thm:congr}]
Let $\rho=(1\,2\,\cdots\,p)\in S_c$ (recall $c\ge p$) and let $X=\{\pi\in S_c:\pi^k=1\}$.

\emph{Step 1: an action of order $p$.} The cyclic group $\langle\rho\rangle$, of order $p$, acts on $S_c$ by conjugation, $\pi\mapsto\rho^i\pi\rho^{-i}$. Conjugation preserves the order and the sign, so this action restricts to $X$ and preserves $\sgn$. Each orbit has size dividing $p$, hence size $1$ or $p$, and all elements of an orbit have the same sign. Therefore
\begin{equation}\label{eq:orbit}
|X|\equiv|X^\rho|\pmod p,\qquad \sum_{\pi\in X}\sgn\pi\equiv\sum_{\pi\in X^\rho}\sgn\pi\pmod p,
\end{equation}
where $X^\rho$ is the set of fixed points of the action, since the orbits of size $p$ contribute multiples of $p$ to both sums.

\emph{Step 2: the fixed points.} An element $\pi$ is fixed if and only if $\pi$ commutes with $\rho$, i.e.\ $\pi$ lies in the centralizer $Z(\rho)$ of $\rho$ in $S_c$. The group $\langle\rho\rangle\times S_{\{p+1,\dots,c\}}$ centralizes $\rho$ and has order $p\,(c-p)!$; on the other hand the conjugacy class of $\rho$ consists of the $\binom{c}{p}(p-1)!$ $p$-cycles of $S_c$, so $|Z(\rho)|=c!/\big(\binom cp(p-1)!\big)=p\,(c-p)!$. Hence
\[
Z(\rho)=\langle\rho\rangle\times S_{\{p+1,\dots,c\}}=\{\rho^i\pi':\ 0\le i\le p-1,\ \pi'\in S_{\{p+1,\dots,c\}}\}.
\]
Since $\rho^i$ and $\pi'$ have disjoint supports, they commute and $(\rho^i\pi')^k=\rho^{ik}\pi'^k$, which is the identity if and only if $\rho^{ik}=1$ and $\pi'^k=1$. Moreover, $\sgn(\rho^i\pi')=\sgn(\rho^i)\sgn(\pi')=\sgn(\pi')$, because $\rho^i$ is a power of an odd-length cycle and hence even. Thus
\[
X^\rho=\{\rho^i\pi':\ p\mid ik,\ \pi'\in S_{\{p+1,\dots,c\}},\ \pi'^k=1\}.
\]

\emph{Step 3: counting.} If $p\nmid k$, then $p\mid ik$ forces $i=0$, so $X^\rho$ is in sign-preserving bijection with the solutions of $x^k=1$ in $S_{c-p}$; by \eqref{eq:orbit}, $r_k(c)\equiv r_k(c-p)$ and $T_k(c)\equiv T_k(c-p)\pmod p$. If $p\mid k$, all $p$ values of $i$ occur, so $|X^\rho|=p\,r_k(c-p)$ and $\sum_{X^\rho}\sgn=p\,T_k(c-p)$, both $\equiv0\pmod p$.

\emph{Step 4: even solutions.} We have $e_k(c)=\tfrac12(r_k(c)+T_k(c))$ for all $c$. Since $p$ is odd, $2$ is invertible modulo $p$, and the statements for $e_k$ follow from those for $r_k$ and $T_k$.
\end{proof}

\begin{proof}[Proof of Corollary~\ref{cor:radical}]
Assume $T_k(c)=0$ with $c\ge2$, and let $q$ be a prime dividing $c(c-1)$. If $q=2$, then $q\mid k$ because $k$ is even. Let $q$ be odd and suppose, for a contradiction, that $q\nmid k$. Since $q$ divides $c$ or $c-1$, we have $c\equiv0$ or $c\equiv1\pmod q$; in the first case $c\ge q$, and in the second case $c\ge q+1$, because $c-1\ge1$ is divisible by $q$. Hence part~(1)(b) of Theorem~\ref{thm:congr} (with $p=q$) can be applied repeatedly to $c,\,c-q,\,c-2q,\dots$ (the last application being to an integer in $\{q,q+1\}$), and we obtain
\[
T_k(c)\equiv T_k(0)=1\quad\text{or}\quad T_k(c)\equiv T_k(1)=1\pmod q,
\]
respectively. This contradicts $T_k(c)=0$. Hence every prime divisor of $c(c-1)$ divides $k$.

Let $P$ be the set of prime divisors of $k$. By what we just proved,
if $T_k(c)=0$ with $c\ge2$, then every prime factor of $c-1$ and of $c$ lies in $P$;
that is, $(c-1,c)$ is a pair of consecutive positive integers whose prime factors all
lie in $P$. By a theorem of St\o rmer \cite{Stormer} (see also \cite[Section~1]{Lehmer}), for every finite set $P$ of primes
there are only finitely many such pairs, and they can be found by solving finitely many
Pell equations. Hence $\{c:T_k(c)=0\}$ is finite (recall that $T_k(0)=T_k(1)=1$). If $k=2^a$, then $P=\{2\}$, so $c-1$ and $c$
are both powers of $2$, which forces $c=2$; conversely, $T_k(2)=0$ by
Theorem~\ref{thm:balance}.
\end{proof}

Since $T_k(c)$ depends only on the set $\{g\le c: g\mid k\}$ (a permutation of $[c]$ satisfies $\pi^k=1$ if and only if each of its cycle lengths, which are at most $c$, divides $k$), the zeros of $T_k$ can be determined by a finite computation for each fixed set of primes. For the example $k=2^a3^b$ given after Conjecture~\ref{conj:T}, the pairs of consecutive positive integers with no prime factors other than $2$ and $3$ are $(1,2)$, $(2,3)$, $(3,4)$ and $(8,9)$. Indeed, two consecutive integers are coprime, so, apart from the pair $(1,2)$, one of them is $2^s$, and the other is $3^t$ with $s,t\ge1$. If $2^s-3^t=1$, then $2^s\equiv1\pmod 3$, so $s$ is even and $(2^{s/2}-1)(2^{s/2}+1)=3^t$; two powers of $3$ that differ by $2$ are $1$ and $3$, so $s=2$ and the pair is $(3,4)$. If $3^t-2^s=1$ and $s\ge2$, then $3^t\equiv1\pmod 4$, so $t$ is even and $(3^{t/2}-1)(3^{t/2}+1)=2^s$; two powers of $2$ that differ by $2$ are $2$ and $4$, so $t=2$ and the pair is $(8,9)$; and if $s=1$, the pair is $(2,3)$. Hence only $c\in\{2,3,4,9\}$ can occur. For $c\le9$ the set $\{g\le c:g\mid k\}$, and hence $T_k(c)$, depends only on $\min(a,3)$ and $\min(b,2)$, since the relevant divisors are $1,2,3,4,6,8,9$. Table~\ref{tab:23} lists the nine resulting cases, and the stated description follows.

\begin{table}[ht]
\centering
\begin{tabular}{@{}lrrrrrrrrr@{}}
\toprule
$k$ & 2 & 4 & 8 & 6 & 12 & 24 & 18 & 36 & 72\\
\midrule
$T_k(3)$ & $-2$ & $-2$ & $-2$ & $0$ & $0$ & $0$ & $0$ & $0$ & $0$\\
$T_k(4)$ & $-2$ & $-8$ & $-8$ & $6$ & $0$ & $0$ & $6$ & $0$ & $0$\\
$T_k(9)$ & $28$ & $6832$ & $-38528$ & $-1764$ & $5040$ & $-40320$ & $38556$ & $45360$ & $0$\\
\bottomrule
\end{tabular}
\caption{The values of $T_k(c)$, $c\in\{3,4,9\}$, for $k=2^a3^b$ with $1\le a\le3$ and $0\le b\le2$. These determine $T_k(c)$ for every $k=2^a3^b$, $a\ge1$, and $c\le 9$.}
\label{tab:23}
\end{table}

\begin{proof}[Proof of Proposition~\ref{prop:conj40}]
Fix $c\ge2$ and let $L_c=\operatorname{lcm}(1,\dots,c)$. An integer $g\le c$ divides $k$ if
and only if it divides $\gcd(k,L_c)$, because $g\mid L_c$. Since $T_k(c)$ depends only on
the set $\{g\le c:g\mid k\}$, we get $T_k(c)=T_{\gcd(k,L_c)}(c)$. Moreover, $c$ and $c-1$
are coprime and both divide $L_c$, so $c(c-1)\mid L_c$, and therefore $c(c-1)\mid k$ if
and only if $c(c-1)\mid\gcd(k,L_c)$. Finally, if $k$ is even, then $\gcd(k,L_c)$ is an
even divisor of $L_c$. Consequently, Conjecture~\ref{conj:T} holds for $c$ if and only if
$c(c-1)\mid d$ for every even divisor $d$ of $L_c$ with $T_d(c)=0$, which is a finite
verification. We carried it out for $2\le c\le40$ in exact integer arithmetic, computing
$T_d(c)$ from the recurrence
\[
T_k(c)=\sum_{\substack{g\mid k\\ 1\le g\le c}}(-1)^{g+1}\frac{(c-1)!}{(c-g)!}\,T_k(c-g),\qquad T_k(0)=1,
\]
which follows by removing the cycle containing $c$ (it has length $g\mid k$, sign
$(-1)^{g+1}$, and can be chosen in $(c-1)!/(c-g)!$ ways). For $c=40$ there are $30720$
even divisors of $L_{40}$ to check. The programs are available as ancillary files of the arXiv version of this paper.
\end{proof}

The same argument shows that the least even $k$ with $T_k(c)=0$ is an even divisor of
$L_c$, since $\gcd(k,L_c)\le k$. These least values are listed in Table~\ref{tab:zeros}; in
particular, $T_k$ has zeros at every $c\le40$ for suitable even $k$. Note that
$L_9=L_{10}=2520$; thus, for $c\le10$, it suffices to consider the even $k\le2520$.

\begin{table}[ht]
\centering
\small\setlength{\tabcolsep}{3pt}
\begin{tabular}{@{}lrrrrrrrrrr@{}}
\toprule
$c$ & 3 & 4 & 5 & 6 & 7 & 8 & 9 & 10 & 11 & 12\\
least $k$ & 6 & 12 & 20 & 30 & 42 & 168 & 72 & 90 & 660 & 660\\
\midrule
$c$ & 13 & 14 & 15 & 16 & 17 & 18 & 19 & 20 & 21 & 22\\
least $k$ & 780 & 5460 & 420 & 720 & 12240 & 21420 & 37620 & 47880 & 2520 & 27720\\
\midrule
$c$ & 23 & 24 & 25 & 26 & 27 & 28 & 29 & 30 & 31 & 32\\
least $k$ & 182160 & 212520 & 92400 & 163800 & 196560 & 83160 & 1607760 & 1741740 & 3723720 & 4062240\\
\midrule
$c$ & 33 & 34 & 35 & 36 & 37 & 38 & 39 & 40 & &\\
least $k$ & 1441440 & 471240 & 942480 & 3112200 & 103030200 & 97435800 & 13693680 & 8288280 & &\\
\bottomrule
\end{tabular}
\caption{The least even $k$ with $T_k(c)=0$, for $3\le c\le40$; in every case $c(c-1)\mid k$. For $c\in\{2,\dots,7,9,10\}$ one has $T_k(c)=0$ if and only if $c(c-1)\mid k$, for every even $k$; for the other values of $c\le40$ the condition $c(c-1)\mid k$ is not sufficient.}
\label{tab:zeros}
\end{table}

\section{Proof of Theorem~\ref{thm:rigid} and Corollary~\ref{cor:An}}\label{sec:rigid}
\begin{proof}[Proof of Theorem~\ref{thm:rigid}]
For each $\ell$, let $X_\ell(\sigma)$ be the set of points lying in $\ell$-cycles of $\sigma$.
By Lemma~\ref{lem:cyclepower}, every $k$th root $\tau$ of $\sigma$ leaves each $X_\ell(\sigma)$
invariant, and $\tau\mapsto(\tau|_{X_\ell(\sigma)})_\ell$ is a bijection between the $k$th roots
of $\sigma$ and the families of $k$th roots of the restrictions $\sigma|_{X_\ell(\sigma)}$, which
have type $(\ell)^{c_\ell(\sigma)}$; moreover, $\sgn\tau=\prod_\ell\sgn(\tau|_{X_\ell(\sigma)})$.
Hence all roots of $\sigma$ have the same parity if and only if, for each $\ell$, all roots of
a permutation of type $(\ell)^{c_\ell}$ have the same parity, i.e.\ $|D_k(\ell,c_\ell)|=r_k(\ell,c_\ell)$
for every $\ell$; and in that case the common sign is $\prod_\ell\varepsilon_\ell$, where
$\varepsilon_\ell$ is the common sign of the roots of type $(\ell)^{c_\ell}$.

\emph{Even $\ell$.} By Lemma~\ref{lem:even-ell}, the condition always holds, with $\varepsilon_\ell=(-1)^{c_\ell/2^a}$.

\emph{Odd $\ell$.} Let $g_0=\dd{\ell}{k}$, $k'=k/g_0$, $L=\ell g_0$ and $j=c_\ell/g_0$ (an integer, as $\sigma$ has roots). By Lemma~\ref{lem:reduction}, $r_k(\ell,c_\ell)=\lambda\,r_{k'}(L,j)$ and $D_k(\ell,c_\ell)=\lambda\,D_{k'}(L,j)$ with the same $\lambda>0$, so the condition is equivalent to $|D_{k'}(L,j)|=r_{k'}(L,j)$. Since $L$ is odd and coprime to $k'$, Lemma~\ref{lem:coprime} gives
\[
r_{k'}(L,j)=\sum_{\pi^{k'}=1}L^{\,j-\cyc(\pi)},\qquad D_{k'}(L,j)=\sum_{\pi^{k'}=1}\sgn(\pi)\,L^{\,j-\cyc(\pi)},
\]
with all weights $L^{j-\cyc(\pi)}$ positive. Hence $|D_{k'}(L,j)|=r_{k'}(L,j)$ if and only if all solutions of $x^{k'}=1$ in $S_j$ have the same sign. Since $g_0$ is odd (its prime divisors divide $\ell$) and $k$ is even, $k'$ is even. If $j\ge2$, the identity and the transposition $(1\,2)$ are solutions of $x^{k'}=1$ in $S_j$ of opposite signs, so the condition fails. If $j\le1$, the only solution is the identity, and the condition holds. Thus, for odd $\ell$, the condition is $j\in\{0,1\}$, i.e.\ $c_\ell\in\{0,g_0\}=\{0,\dd{\ell}{k}\}$.

\emph{The sign.} Suppose the condition holds and let $\ell$ be odd. If $c_\ell=0$, there is nothing to restrict. If $c_\ell=g_0$, then $j=1$ and, by Lemma~\ref{lem:structure} and Lemma~\ref{lem:cyclepower}, every root restricted to $X_\ell(\sigma)$ is a single cycle of length $g_0\ell$ (no other $g\in G_k(\ell)$ satisfies $g\le c_\ell=g_0$). This length is odd, so $\varepsilon_\ell=+1$. Therefore the common sign of all roots is $\prod_{\ell\ \mathrm{even}}(-1)^{c_\ell/2^a}$.
\end{proof}

\begin{example}
Let $k=2$. For $\sigma=(1\,2)(3\,4)(5\,6\,7)$ we have $c_2=2$, $c_3=1=\dd{3}{2}$, so all square roots of $\sigma$ have sign $(-1)^{2/2}=-1$; indeed they are $(1\,3\,2\,4)(5\,7\,6)$ and $(1\,4\,2\,3)(5\,7\,6)$. For $\sigma=(1\,2)(3\,4)(5\,6\,7)(8\,9\,10)$, instead, $c_3=2$ and $\sigma$ has square roots of both parities, e.g.\ $(1\,3\,2\,4)(5\,7\,6)(8\,10\,9)$ (odd) and $(1\,3\,2\,4)(5\,8\,6\,9\,7\,10)$ (even).
\end{example}

\begin{proof}[Proof of Corollary~\ref{cor:An}]
Since $k$ is even, $\sgn(\tau^k)=1$ for every $\tau$, so every $k$th power is even. By the criterion of Knopfmacher and Warlimont (see \cite[Theorem~4.43, p.~158]{Wilf}), $\sigma$ has $k$th roots if and only if $\dd{\ell}{k}\mid c_\ell(\sigma)$ for every $\ell$. Assume this holds. Then $\sigma$ fails to have an even $k$th root if and only if all its $k$th roots are odd, that is, if and only if all its roots have the same parity and the common sign is $-1$. By Theorem~\ref{thm:rigid}, this happens if and only if $c_\ell(\sigma)\in\{0,\dd{\ell}{k}\}$ for every odd $\ell$ and $\sum_{\ell\ \mathrm{even}}c_\ell(\sigma)/2^a$ is odd. Since $\dd{\ell}{k}\mid c_\ell(\sigma)$, the condition $c_\ell(\sigma)\in\{0,\dd{\ell}{k}\}$ is equivalent to $c_\ell(\sigma)\le\dd{\ell}{k}$. Negating, $\sigma$ has an even $k$th root if and only if $c_\ell(\sigma)>\dd{\ell}{k}$ for some odd $\ell$ or $\sum_{\ell\ \mathrm{even}}c_\ell(\sigma)/2^a$ is even.
\end{proof}

\subsection{\texorpdfstring{Permutations with a unique $k$th root}{Permutations with a unique kth root}}
The reduction of Section~\ref{sec:prelim} also characterizes the permutations with a
unique $k$th root. For $k$ even, the least prime divisor of $k$ is $2$, and the condition
below says that all cycle lengths of $\sigma$ are odd, coprime to $k$ and pairwise
distinct; in particular $c_\ell(\sigma)\le1=\dd{\ell}{k}$ for every odd $\ell$ and
$\sigma$ has no cycles of even length, so these permutations are an extreme case of
Theorem~\ref{thm:rigid}.

\begin{corollary}\label{cor:unique}
Let $k\ge2$ and let $p$ be the least prime divisor of $k$. A permutation $\sigma$ has exactly one $k$th root if and only if, for every $\ell$ with $c_\ell(\sigma)>0$, we have $\gcd(\ell,k)=1$ and $c_\ell(\sigma)<p$. In that case the unique root is $\sigma^{k^*}$, where $kk^*\equiv1\pmod{\operatorname{ord}\sigma}$; in particular, it is a power of $\sigma$. The exponential generating function of the number $u_k(n)$ of permutations of $[n]$ with a unique $k$th root is
\[
\sum_{n\ge0}u_k(n)\frac{x^n}{n!}=\prod_{\gcd(\ell,k)=1}\ \sum_{i=0}^{p-1}\frac{1}{i!}\Big(\frac{x^\ell}{\ell}\Big)^{i}.
\]
\end{corollary}

\begin{proof}
By \eqref{eq:product}, $r_k(\sigma)=1$ if and only if $r_k(\ell,c)=1$ for every $\ell$
with $c=c_\ell(\sigma)>0$. Fix such an $\ell$, let $g_0=\dd{\ell}{k}$, and assume
$g_0\mid c$ (otherwise $r_k(\ell,c)=0$ by Lemma~\ref{lem:reduction}). By
Lemma~\ref{lem:cyclepower} and \cite[Theorem~1]{LMR}, every family
$(s_g)_{g\in G_k(\ell)}$ of non-negative integers with $\sum_g g\,s_g=c$ is the cycle
type $\prod_g(g\ell)^{s_g}$ of at least one $k$th root of a permutation of type
$(\ell)^c$, and every such root arises in this way. By Lemma~\ref{lem:structure}, every
$g\in G_k(\ell)$ is a multiple of $g_0$. If some $g\in G_k(\ell)$ satisfies $g_0<g\le c$,
the families with $s_{g_0}=c/g_0$, and with $s_g=1$ and $s_{g_0}=(c-g)/g_0$, give two
distinct roots. Otherwise, all roots have type $(g_0\ell)^{c/g_0}$, and by
\eqref{eq:onetype} their number is the number of permutations of $[c]$ of type
$(g_0)^{c/g_0}$ times $\ell^{(g_0-1)c/g_0}$; this is $1$ if $g_0=1$, and at least
$\ell\ge2$ if $g_0\ge2$ (because $\dd{1}{k}=1$). Since $g_0=1$ means $\gcd(\ell,k)=1$,
and in that case $G_k(\ell)$ is the set of divisors of $k$ (Lemma~\ref{lem:structure}),
whose second smallest element is $p$, we conclude that $r_k(\ell,c)=1$ if and only if
$\gcd(\ell,k)=1$ and $c<p$.

If these conditions hold for every $\ell$ with $c_\ell(\sigma)>0$, then
$\operatorname{ord}\sigma=\operatorname{lcm}\{\ell:c_\ell(\sigma)>0\}$ is coprime to $k$,
and $(\sigma^{k^*})^k=\sigma$; so $\sigma^{k^*}$ is the unique root. Finally,
Lemma~\ref{lem:expformula} with $w_\ell=1$, $A_\ell=\{0,\dots,p-1\}$ if $\gcd(\ell,k)=1$
and $A_\ell=\{0\}$ otherwise, gives the generating function.
\end{proof}

For $k=2$ we recover the classical fact that a permutation has a unique square root if and only if its cycle lengths are odd and pairwise distinct (see the comments to \cite[A088994]{OEIS}); the generating function is $\prod_{\ell\ \mathrm{odd}}(1+x^\ell/\ell)$, and the counting sequence begins
\[1,1,0,2,8,24,144,720,8448,64512,576000,\dots\]

\begin{example}
For $k=6$ (so $p=2$) the condition is that all cycle lengths are coprime to $6$ and pairwise distinct; for $k=9$ (so $p=3$) all cycle lengths must be coprime to $3$, each occurring at most twice. For instance $(1\,2)(3\,4)$ has the unique $9$th root $(1\,2)(3\,4)$ itself, while $(1\,2)(3\,4)(5\,6)$ has the $9$th roots $(1\,2)(3\,4)(5\,6)$ and the eight $6$-cycles whose $9$th power is $(1\,2)(3\,4)(5\,6)$.
\end{example}

\section{\texorpdfstring{Proofs of Theorems~\ref{thm:Aodd} and \ref{thm:A2r}: powers in the alternating group}{Proofs of Theorems 1.8 and 1.9: powers in the alternating group}}\label{sec:An}

Throughout this section $\POW_m(n)$ denotes the set of $m$th powers in $S_n$, so that $p_m(n)=|\POW_m(n)|/n!$, and for $n\ge2$ we have $|A_n|=n!/2$, so $p^A_m(n)=2|A^m_n|/n!$. We use that, for $m=q^r$ with $q$ prime,
\begin{equation}\label{eq:pm}
p_m(qj)=p_m(qj+1)=\cdots=p_m(qj+q-1)\qquad(j\ge0)
\end{equation}
\cite{Chernoff94,LMR}. We also use the following special case of the criterion for $m$th roots: if $m=q^r$ with $q$ prime, then $\dd{\ell}{m}=m$ if $q\mid\ell$ and $\dd{\ell}{m}=1$ otherwise; hence
\begin{equation}\label{eq:POWq}
\sigma\in\POW_m(n)\iff m\mid c_\ell(\sigma)\ \text{for every $\ell$ divisible by $q$}.
\end{equation}

\begin{proof}[Proof of Theorem~\ref{thm:Aodd}]
\emph{Step 1: $A_n^m$ for $m$ odd.} If $\tau\in A_n$, then $\tau^m\in A_n$ and $\tau^m\in\POW_m(n)$. Conversely, if $\sigma=\tau^m$ is even with $\tau\in S_n$, then $\sgn\tau=(\sgn\tau)^m=\sgn\sigma=1$, so $\tau\in A_n$. Hence $A_n^m=\POW_m(n)\cap A_n$. Let
\[
s(n)=\sum_{\sigma\in\POW_m(n)}\sgn\sigma=|\POW_m(n)\cap A_n|-|\POW_m(n)\setminus A_n|.
\]
Then $|A^m_n|=\tfrac12(|\POW_m(n)|+s(n))$, and therefore
\begin{equation}\label{eq:pAodd}
p^A_m(n)=\frac{2|A_n^m|}{n!}=p_m(n)+\frac{s(n)}{n!}\qquad(n\ge2).
\end{equation}

\emph{Step 2: the generating function of $s(n)$.} By \eqref{eq:POWq} (with $q=p$), $\POW_m(n)$ consists of the permutations with $c_\ell\in m\N_0=\{0,m,2m,\dots\}$ when $p\mid\ell$ and $c_\ell$ arbitrary when $p\nmid\ell$. The sign of $\sigma$ is $\prod_C(-1)^{|C|+1}$. Lemma~\ref{lem:expformula} with $w_\ell=(-1)^{\ell+1}$ gives
\[
\sum_{n\ge0}s(n)\frac{x^n}{n!}=\prod_{p\nmid\ell}\exp\Big(\frac{(-1)^{\ell+1}x^\ell}{\ell}\Big)\ \prod_{i\ge1}\exp_m\Big(\frac{(-1)^{pi+1}x^{pi}}{pi}\Big),
\]
since $\sum_{t\in m\N_0}y^t/t!=\exp_m(y)$. As $p$ is odd, $(-1)^{pi+1}=(-1)^{i+1}$, so the second product is the one appearing in the definition of $\kappa_j$.

\emph{Step 3: the first product.} Using $\sum_{\ell\ge1}(-1)^{\ell+1}x^\ell/\ell=\log(1+x)$ and again $(-1)^{pi+1}=(-1)^{i+1}$,
\[
\sum_{p\nmid\ell}\frac{(-1)^{\ell+1}x^\ell}{\ell}=\log(1+x)-\sum_{i\ge1}\frac{(-1)^{i+1}x^{pi}}{pi}=\log(1+x)-\frac1p\log(1+x^p).
\]
Hence the first product equals $(1+x)\exp\big(-\tfrac1p\log(1+x^p)\big)=(1+x)(1+x^p)^{-1/p}$, and
\[
\sum_{n\ge0}s(n)\frac{x^n}{n!}=(1+x)\sum_{j\ge0}\kappa_jx^{pj}=\sum_{j\ge0}\kappa_j\big(x^{pj}+x^{pj+1}\big).
\]
Because $p\ge3$, the exponents $pj$ and $pj+1$ ($j\ge0$) are pairwise distinct and the exponents $pj+i$ with $2\le i\le p-1$ do not occur. Thus $s(n)/n!=\kappa_j$ if $n\in\{pj,pj+1\}$ and $s(n)=0$ if $n=pj+i$ with $2\le i\le p-1$.

\emph{Step 4: conclusion.} By \eqref{eq:pm}, $p_m(pj+i)=p_m(pj)$ for $0\le i\le p-1$. Substituting this and Step~3 into \eqref{eq:pAodd} gives the theorem.
\end{proof}

In particular, the signed number of $m$th powers in $S_n$ vanishes unless $n\equiv0,1\pmod p$.

\begin{example}
Let $m=p=3$. Then $\kappa_0=1$, and for $n=3$ Theorem~\ref{thm:Aodd} gives $p^A_3(3)=p_3(3)+\kappa_1$. The cubes in $S_3$ are the identity and the three transpositions, so $p_3(3)=2/3$, while $A_3^3=\{1\}$, so $p^A_3(3)=1/3$ and $\kappa_1=-1/3$. Directly from the definition, $[x^3]\big((1+x^3)^{-1/3}\exp_3(x^3/3)\big)=-1/3$, in agreement.
\end{example}

\begin{proof}[Proof of Theorem~\ref{thm:A2r}]
\emph{Step 1: the set $B_n$.} Since $m$ is even, every $m$th power is even, so $\POW_m(n)\subseteq A_n$, and $\sigma\in A_n^m$ if and only if $\sigma$ has an even $m$th root. Hence $|A_n^m|=|\POW_m(n)|-|B_n|$, where $B_n$ is the set of $m$th powers all of whose $m$th roots are odd. Therefore, for $n\ge2$,
\begin{equation}\label{eq:pA2}
p^A_m(n)=\frac{2|A^m_n|}{n!}=2p_m(n)-2\,\frac{|B_n|}{n!}.
\end{equation}

\emph{Step 2: description of $B_n$.} Here $a=\nu_2(m)=r$ and $2^a=m$. By \eqref{eq:POWq} (with $q=2$), $\dd{\ell}{m}=1$ for odd $\ell$ and $\dd{\ell}{m}=m$ for even $\ell$. By (the proof of) Corollary~\ref{cor:An}, $\sigma\in B_n$ if and only if
\[
c_\ell\le1\ \text{for odd $\ell$},\qquad m\mid c_\ell\ \text{for even $\ell$},\qquad \sum_{\ell\ \mathrm{even}}c_\ell/m\ \text{is odd}.
\]

\emph{Step 3: generating function.} Fix $\epsilon\in\{1,-1\}$ and apply Lemma~\ref{lem:expformula} (with variable $x$) to the sequences
\[
a_\ell(i)=\begin{cases}1,& i\in\{0,1\},\\ 0,& i\ge2,\end{cases}\quad(\ell\ \text{odd}),\qquad
a_\ell(i)=\begin{cases}\epsilon^{\,i/m},& m\mid i,\\ 0,& m\nmid i,\end{cases}\quad(\ell\ \text{even}).
\]
Then $\prod_\ell a_\ell(c_\ell(\pi))$ vanishes unless $\pi$ satisfies the first two conditions of Step~2, in which case it equals $\epsilon^{\sum_{\ell\ \mathrm{even}}c_\ell(\pi)/m}$. On the right-hand side of Lemma~\ref{lem:expformula}, each odd $\ell$ contributes the factor $1+x^\ell/\ell$, and each even $\ell$ contributes $\sum_{t\ge0}\epsilon^{t}(x^\ell/\ell)^{tm}/(tm)!$, which is $\exp_m(x^\ell/\ell)$ if $\epsilon=1$ and $\widetilde{\exp}_m(x^\ell/\ell)$ if $\epsilon=-1$. Hence
\[
\sum_{n\ge0}\frac{x^n}{n!}\sum_{\pi}\epsilon^{\sum_{\ell\ \mathrm{even}}c_\ell(\pi)/m}=
\begin{cases}O(x)\,E(x),&\epsilon=1,\\ O(x)\,\widetilde E(x),&\epsilon=-1,\end{cases}
\]
where the inner sum runs over the permutations $\pi\in S_n$ satisfying the first two conditions of Step~2. Taking half the difference of the two cases counts exactly those $\pi$ for which, in addition, $\sum_{\ell\ \mathrm{even}}c_\ell(\pi)/m$ is odd, that is, the elements of $B_n$. Therefore
\[
B(x):=\sum_{n\ge0}|B_n|\frac{x^n}{n!}=O(x)\,\tfrac12\big(E(x)-\widetilde E(x)\big).
\]
Together with \eqref{eq:pA2} this proves the formula for $p^A_m(n)$.

\emph{Step 4: the coefficients.} Write $b_n=[x^n]B(x)$, $O(x)=\sum_nO_nx^n$ and $\tfrac12(E-\widetilde E)=\sum_ne_nx^n$. Only even cycle lengths occur in $E$ and $\widetilde E$, so $e_n=0$ for $n$ odd. The smallest $n$ for which some configuration of even cycles has $\sum c_\ell/m$ odd is $n=2m$, attained only by $c_2=m$ (indeed, if $\sum_{\ell\ \mathrm{even}}c_\ell/m$ is odd then $\sum_{\ell\ \mathrm{even}}c_\ell\ge m$, so $n=\sum_{\ell\ \mathrm{even}}\ell c_\ell\ge2\sum_{\ell\ \mathrm{even}}c_\ell\ge2m$, with equality only if $c_2=m$ and $c_\ell=0$ for even $\ell\ge4$); hence $e_{2i}=0$ for $i<m$ and
\[
e_{2m}=\frac{1}{m!}\Big(\frac{x^2}{2}\Big)^m\Big/x^{2m}=\frac{1}{2^m\,m!}.
\]
Expanding the product $(1+x)(1+x^3/3)(1+x^5/5)\cdots$ we get
\[
O(x)=1+x+\tfrac13x^3+\tfrac13x^4+\cdots,\qquad\text{so}\quad O_0=O_1=1,\ O_2=0,\ O_3=\tfrac13.
\]
Since $e_n=0$ for odd $n$,
\[
b_{2j+1}-b_{2j}=\sum_{i\ge0}O_{2i+1}e_{2j-2i}-\sum_{i\ge0}O_{2i}e_{2j-2i}=\sum_{i\ge0}(O_{2i+1}-O_{2i})\,e_{2j-2i}.
\]

\emph{Step 5: conclusion.} By \eqref{eq:pm} with $q=2$, $p_m(2j)=p_m(2j+1)$, so by \eqref{eq:pA2}
\[
p^A_m(2j)-p^A_m(2j+1)=2\,(b_{2j+1}-b_{2j})\qquad(j\ge1).
\]
If $j\le m$, then $e_{2j-2i}\ne0$ requires $j-i\ge m$, i.e.\ $i\le j-m\le0$; so only $i=0$ may contribute, and its contribution is $(O_1-O_0)e_{2j}=0$. Hence $p^A_m(2j)=p^A_m(2j+1)$ for $1\le j\le m$. If $j=m+1$, only $i\in\{0,1\}$ may contribute; the term $i=0$ vanishes as before, and the term $i=1$ equals $(O_3-O_2)e_{2m}=e_{2m}/3$. Therefore
\[
p^A_m(2m+2)-p^A_m(2m+3)=\frac{2e_{2m}}{3}=\frac{2}{3\cdot2^m m!}.\qedhere
\]
\end{proof}

For $m=2$ we have $\exp_2=\cosh$ and $\widetilde{\exp}_2=\cos$, and Theorem~\ref{thm:A2r} recovers Pournaki's generating function \cite{Pournaki}, and the difference $p^A_2(6)-p^A_2(7)=1/12$ agrees with the value $2/(3\cdot2^mm!)$ for $m=2$. While Theorem~\ref{thm:A2r} guarantees $p^A_m(2j)=p^A_m(2j+1)$ for $1\le j\le m$, numerically, for $m\in\{2,4,8\}$ and $2j+1\le21$, the equality fails for every $j>m$.

\section*{Acknowledgements}
The author used the large language model Claude (Anthropic) to: write and run computer
programs that check the numerical results (Tables~\ref{tab:23} and~\ref{tab:zeros}, the values of $p^A_m(n)$ and
the sequence $u_2(n)$), and that carry out the finite verification in the proof of Proposition~\ref{prop:conj40}; expand the proofs with additional intermediate steps and examples;
prepare the figures; and edit the
English and the \LaTeX{} source. The author checked every step of the proofs and every citation, and takes
full responsibility for the content.

\end{document}